\documentclass[preprint,authoryear]{elsarticle}

\usepackage[english]{babel}			
\usepackage{amsmath,amsthm,amsfonts,amssymb}
\usepackage{paralist}
\usepackage{pgfplots}
\usepackage{tikz}
\usepackage[top=2.5cm, left=2.5cm, right=2.5cm, bottom=2cm]{geometry}

\def\N{\mathbb{N}}
\def\P{\mathbb{P}}
\def\R{\mathbb{R}}

\def\Po{\text{Po}}
\def\diam{\text{diam}}
\def\Pol{\text{Pole}}
\def\NA{\text{NA}}

\newtheorem{env_theo}{Theorem}[section]
\newtheorem{env_lemma}[env_theo]{Lemma}
\newtheorem{env_cor}[env_theo]{Corollary}
\newtheorem{env_prop}[env_theo]{Proposition}

\begin{document}

\begin{frontmatter}

\title{The limit distribution of the largest interpoint distance for distributions supported by an ellipse and generalizations}

\author[kit]{Michael Schrempp}
\ead{schrempp@kit.edu}

\address[kit]{Karlsruhe Institute of Technology, Institute of Stochastics \\
Englerstr. 2, D-76131 Karlsruhe, Germany}

\begin{abstract}
We study the asymptotic behaviour of the maximum interpoint distance of random  points in a planar bounded set with an unique major axis and a boundary behaving like an ellipse at the endpoints. Our main result covers the case of uniformly distributed points in an ellipse.
\end{abstract}

\end{frontmatter}


\section{Introduction}
For some fixed integer $d \ge 2$, let $X_1,X_2,\ldots$ be a sequence of independent and identically distributed (i.i.d.) $d$-dimensional random vectors, defined  on a common probability space $(\Omega,{\cal A},\P)$.
Writing $|\cdot |$ for the Euclidean norm on $\R^d$, the  convergence in distribution of the suitably normalized maximum interpoint distance
$$
M_n := \max\limits_{1 \le i < j \le n}|X_i - X_j|
$$
has been a topic of interest for more than 20 years. Results obtained so far are mostly for the case that the distribution $\P_{X_1}$ of $X_1$ is spherically symmetric, and they may roughly be classified according to whether $\P_{X_1}$ has an unbounded or a bounded support. If $X_1$ has a spherically symmetric normal distribution, \cite{Matthews1993} obtained a Gumbel limit distribution
for $M_n$. \cite{Henze1996} generalized this result to the case that $X_1$ has a spherically symmetric Kotz distribution.
An even more general spherically symmetric setting has recently  been studied by \cite{RaoJanson2015}. In the unbounded case, \cite{HenzeLao2010} obtained a (non-Gumbel) limit distribution of $M_n$ if the distribution of $X_1$ is power-tailed spherically decomposable. This case covers certain long-tailed spherically symmetric distributions for $X_1$.
Finally, \cite{Demichel2014} proved several results for the diameter of an elliptical cloud.

If $\P_{X_1}$ has a bounded support, \cite{Appel2002} obtained a convolution of two Weibull distributions as limit law of $M_n$
if $X_1$ has uniform distribution in a planar set with
unique major axis and sub-$\sqrt{x}$ decay of its boundary at the endpoints. Moreover, they derived bounds for the limit law of $M_n$ if $X_1$ has a uniform distribution in an ellipse. \cite{Lao2010}, and \cite{MayerMol2007} obtained Weibull limit distributions for $M_n$  under very general settings if the distribution of $X_1$ is supported by the $d$-dimensional unit ball $\mathbb{B}^d$ for $d \ge 2$ (including the case of a uniform distribution). \cite{Lao2010} obtained limit laws for $M_n$ if $\P_{X_1}$ is uniform or non-uniform in the unit square, uniform in regular polygons, or uniform in the unit $d$-cube, $d \ge 2$.
Moreover, if $\P_{X_1}$ is uniform in a proper ellipse, she improved the lower bound on the limit distribution of $M_n$ given in \cite{Appel2002}. The limit behaviour of $M_n$ if $\P_{X_1}$ is uniform in a proper ellipse has been an open problem for many years. Without giving a proof, \cite{RaoJanson2015} state that $n^{2/3}(2-M_n)$ has a limit distribution (involving two independent Poisson processes) if $X_1$ has a uniform distribution in a proper ellipse with major axis $2$. We generalize this result to the case that the distribution is uniform or non-uniform over a planar bounded set satisfying certain regularity conditions. Furthermore, the limit distribution of $M_n$ will be given in a different way.

In what follows, let $d=2$ and write $\lambda^2$ for Lebesgue measure in the plane. Throughout this work we consider distributions $\P_{X_1}$ with a $\lambda^2$-density $f$ and
compact support $E \subset \R^2$. By our main assumptions (see A1 to A7 in Section~\ref{sec_Assumptions}), $E$ has a unique major axis, and
its boundary decays as fast as $\sqrt{x}$ at the endpoints. In addition, the density $f$ is continuous and bounded away from $0$ at the endpoints. Since the boundary of the unit disk $\mathbb{B}^2$ also decays as fast as $\sqrt{x}$ at the points $(1,0)$ and $(-1,0)$, but $\mathbb{B}^2$ has no unique major axis, this paper can be interpreted as a missing link between the results of \cite{Appel2002} and \cite{Lao2010} resp. \cite{MayerMol2007}.

We also consider a related setting in which the number of points follows a Poisson distribution. To this end, let $\Phi_n$ be a Poisson process in $\R^2$ with intensity measure $n\P_{X_1}$.
The diameter of the support of $\Phi_n$ is denoted by $\diam(\Phi_n)$. With a few assumptions and $\diam(E) := \sup\left\{ |x-y|:x,y \in E\right\}$,
it follows that $\diam(\Phi_n)$ converges almost surely to $\diam(E)$ as $n \to \infty$.
We will show that the assumptions on $E$ and $f$ stated in Section~\ref{sec_Assumptions} imply the weak convergence of $n^{2/3}\big(\diam( E) - \diam(\Phi_n)\big)$.
By use of the De-Poissonization technique we then obtain the weak convergence of
$n^{2/3}\big(\diam(E) - M_n\big)$
towards the same limit distribution.

The rest of the paper is organized as follows: Section~\ref{sec_Assumptions} contains the assumptions on $E$ and $f$, and Section~\ref{sec_Main_Results} includes the main result, which is
Theorem~\ref{satz_Vert_Konv_Hauptsatz}. As a corollary, we obtain the limit distribution of $M_n$ if $X_1$ has a uniform distribution in a proper ellipse. Section~\ref{sec_proofs}
is devoted to proofs. In the final Section~\ref{sec_General_Open_Problems} we indicate possible generalizations and state some open problems.

In the sequel, $A^\circ $ and $\partial A$ stand for the interior and the boundary of a set $A$, respectively, and U$(I)$ denotes the uniform distribution in the interval $I$.
Each unspecified limit refers to $n \to \infty$. Convergence in probability, convergence in distribution and equality in distribution will be
denoted by $\overset{\P}{\longrightarrow }$, $\overset{\mathcal{D}}{\longrightarrow }$ and $\overset{\mathcal{D}}{=}$, respectively.
By $A_n = o_{\P}(1) $ we mean $A_n \overset{\P}{\longrightarrow } 0$. Finally, $f(x) \sim g(x)$ as $x \to x^{*}$ stands for $f(x)/g(x) \rightarrow 1$ as $x \to x^{*}$.


\section{Assumptions and preliminaries}
\label{sec_Assumptions}
We first state the basic assumptions on the set $E$ that supports the distribution of $X_1$. As a bit of notation, we write $B(h) := \{(x,y) \in \R^2: |(x,y)| \le a-h\}$ for the closed circle
centered at the origin with radius $a-h$, where $0 \le h < a$.

\begin{enumerate}
\item[A1)] There is a constant $a>0$ with $\diam(E) = 2a$, and $E$ is oriented in the plane so that $\inf \{ x: (x,y) \in E\} = -a$ and $\sup \{ x: (x,y) \in E\} = a$.
\item[A2)] Putting $U(x,\varepsilon) := \{(r,s) \in \R^2: |(r,s)-(x,0)| < \varepsilon\}$, we have
$\diam\big(E\setminus(U(-a,\varepsilon)\cup U(a,\varepsilon))\big) < 2a$ for each $\varepsilon >0$.
\item[A3)] Writing $Q_i$ for the $i$-th open quadrant in $\R^2$, $i \in \{1,2,3,4\}$, where $Q_1 = \{(x,y): x >0, y>0\}$ and the numbering is anti-clockwise, and putting $E_i := E\cap Q_i$, $i \in \left\{ 1,2,3,4 \right\}$, we have for $i \in \{1,4\}$ and $j \in \{2,3\}$ and each $\varepsilon >0$: $\lambda^2\big( E_i \cap U(a,\varepsilon)\big) > 0$ and $\lambda^2\big( E_j \cap U(-a,\varepsilon)\big) > 0.$
\item[A4)] For some $\nu \in [0,a)$ and continuous functions $g_1:[\nu,a] \rightarrow \R_{\ge 0}$ and $g_4:[\nu,a] \rightarrow \R_{\le 0}$ satisfying $g_1(a) = g_4(a) = 0$, we have $ E^\circ  \cap \{(x,y) \in \R^2 \! : \! x > \nu\} = \{ (x,y) \in \R^2 \! : \! \nu < x < a \text{ and } g_4(x) < y < g_1(x)\}.$ Likewise, for continuous functions $g_2:[-a,- \nu] \rightarrow \R_{\ge 0}, g_3:[-a,- \nu] \rightarrow \R_{\le 0}$ with $g_2(-a) = g_3(-a) = 0$, we have $E^\circ   \cap \{(x,y) \in \R^2 \! : \!  x < -\nu\} = \{ (x,y) \in \R^2: -a < x < -\nu \text{ and } g_3(x) < y < g_2(x)\}.$
\item[A5)] Writing $f_a(x) := \sqrt{a^2-x^2}/2$ for the 'upper boundary function' of an ellipse with major axis $2a$ and minor axis $a$, we assume that for  constants $q_1,q_2,q_3,q_4$ satisfying
\begin{equation} \label{conditiononq}
0 < q_i < 2 , \qquad i \in \left\{ 1,2,3,4 \right\},
\end{equation}
$$
  \frac{g_1(x)}{f_a(x)} \rightarrow  q_1, \ \frac{-g_4(x)}{f_a(x)} \rightarrow  q_4  \quad \text{as } x \rightarrow a, \quad
  \frac{g_2(x)}{f_a(x)} \rightarrow  q_2, \
  \frac{-g_3(x)}{f_a(x)} \rightarrow  q_3  \quad \text{as } x \rightarrow -a.  \
$$
\item[A6)]  For $i \in \left\{ 1,2,3,4\right\}$ and sufficiently small $h$, we have $E_i^\circ \cap \big\{ (x,y) \in \R^2: |x| < \nu \big\} \subset B(h)$. Moreover, $g_i$ has only one point of intersection with $\partial B(h)$. The abscissa of this point is denoted by $\overline x_i(h)$.
\item[A7)] For sufficiently small $\varepsilon$, the density $f$ is continuous on $E\cap U(a,\varepsilon)$ and $E\cap U(-a,\varepsilon)$ and we have $p_1 := p_4 := f\big((a,0)\big) > 0$ and $p_2 := p_3 := f\big((-a,0)\big) > 0$.
\end{enumerate}

Assumption A1 entails no loss of generality since the problem is invariant under rigid motions. A2 means that $(-a,0)$ and $(a,0)$ (henceforth called the
'poles') are the endpoints of the unique major axis of $E$. By A3, the area near the poles is positive in each quadrant.
Assumption A4 means that $E^\circ $ is vertically convex near the poles. In the sequel, $g_1,\ldots,g_4$ will be called the
'boundary functions' of $E$. By A5, these functions decay as fast as $\sqrt{a \pm x}$ at the poles.
For example, the condition in the first quadrant is equivalent to
$
g_1(x)/\sqrt{a-x} \rightarrow   q_1\sqrt{a}/\sqrt{2} > 0$ as $x \rightarrow a$,
which means that $g_1$ actually decays like a square root.
The reason why A5  is formulated in terms of $f_a(x)$ instead of $\sqrt{a\pm x}$ is to facilitate many calculations in Subsection~\ref{subsec_Geom_Ueberl}.
The choice of the factor $1/2$ in $f_a(x)$ is arbitrary but necessary (in fact, it can be any number in the open interval $(0,1)$) in order to have points of intersection of $f_a$ and $\partial B(h)$ in the proofs of Lemma~\ref{lemma_asymp_Verh_Endkappe_EllAehnlich} and Lemma~\ref{lem_Verhalten_gamma_h_mit_Streckung_q_i}.
Notice that (\ref{conditiononq}) ensures the existence of at least one point of intersection  of the boundary functions $g_i$ and the boundary of $B(h)$ for small $h>0$.
If $q_i = 2$ the set $E$ would behave like a circle at the pole in the $i$-th quadrant. This case is explicitly excluded in this work. If $E$ is a circle, we have $q_1 = \ldots = q_4 = 2$. For a circle, the limit distribution of $M_n$ is well-known, see \cite{Lao2010} and \cite{MayerMol2007}. Because of A6, the set $\big\{ z \in E_i^\circ : |z| > a- h\big\}$
consists only of points lying close to a pole for sufficiently small $h$. The notations $p_1 = p_4$ and $p_2 = p_3$ in A7 are redundant but useful, since we hereby avoid a distinction of several cases.


\section{Main results}
\label{sec_Main_Results}
To state the main result, put
\begin{equation}
\label{def_der_Konstanten}
c_i := \frac{2q_i \sqrt{2a}}{3\sqrt{4-q_i^2}}, \quad
\sigma_i := (p_ic_i)^{-2/3}, \quad \tau_i := \frac{3}{2a}c_i, \qquad i \in \left\{ 1,2,3,4\right\}.
\end{equation}
Let $Y_1,Y_2, \ldots $ and $U_1,U_2,\ldots $ be independent random variables, where $Y_1,Y_2,\ldots $ are i.i.d. with a unit exponential distribution, and  $U_1,U_2,\ldots $ are i.i.d. with the uniform
distribution U$([0,1])$. For $m \ge 1$, set $S_m := Y_1+ \ldots + Y_m$. Let $\sigma,\tau >0$, and put
$$
Z_{1,m} := \sigma S_m^{2/3}, \quad Z_{2,m} := U_m \tau \sigma^{1/2} S_m^{1/3} \ (= U_m \tau Z_{1,m}^{1/2}).
$$
The sequence $\mathbf{Z} := (Z_{1,1},Z_{2,1},Z_{1,2},Z_{2,2}, \ldots, Z_{1,m},Z_{2,m}, \ldots )$ defines a $\R^\N$-valued random element, the distribution of which will be denoted by
NA$_\infty(\sigma,\tau)$. Here, the coining $\NA$ stands for 'norm-angle distribution'. Notice that for each $m$ the conditional distribution of $Z_{2,m}$ given $Z_{1,m}$ is uniform on $[0,\tau Z_{1,m}^{1/2}]$. In what follows, we
will write $\mathbf{Z} =: (Z_{1,k},Z_{2,k})_{k \ge 1}$.
Our main result is as follows.

\begin{env_theo}
\label{satz_Vert_Konv_Hauptsatz}
Let $(Z_{1,k}^{i},Z_{2,k}^{i})_{k \ge 1} \sim  \NA_{\infty}(\sigma_i,\tau_i), i \in \{1,2,3,4\},$
be independent random elements of $\R^{\N}$, and put
\begin{eqnarray*}
  S^{i,j} & := &  \min_{ k,l \in \N}\left\{Z_{1,k}^{i} + Z_{1,l}^{j} +\frac{a}{4}\left(Z_{2,k}^{i} - Z_{2,l}^{j}\right)^2 \right\}, \quad (i,j) \in \big\{ (1,3),(2,4)\big\},\\
  S^{i,j} & := &  \min_{ k,l \in \N}\left\{Z_{1,k}^{i} + Z_{1,l}^{j} +\frac{a}{4}\left(Z_{2,k}^{i} + Z_{2,l}^{j}\right)^2 \right\}, \quad (i,j) \in \big\{ (1,2),(3,4)\big\}.
\end{eqnarray*}
Then, under the assumptions A1 - A7, we have
\begin{equation}
\label{eq_Vert_Konv_Aussage}
n^{2/3}\big( 2a - \diam(\Phi_n)\big)
\overset{\mathcal{D}}{\longrightarrow } \min\left\{S^{1,2},S^{1,3},S^{2,4},S^{3,4}\right\}.
\end{equation}
\end{env_theo}
The proof of Theorem~\ref{satz_Vert_Konv_Hauptsatz} is given in Section~\ref{sec_proofs}.
By use of a De-Poissonization theorem by \cite{MayerMol2007}, we can restate this result for the maximum interpoint distance $M_n$ of independent and identically distributed random points.
\begin{env_theo}
\label{satz_Vert_Konv_M_n}
Under A1 to A7 we have $n^{2/3}\left( 2a - M_n\right) \overset{\mathcal{D}}{\longrightarrow } \min\left\{S^{1,2},S^{1,3},S^{2,4},S^{3,4}\right\}$.
\end{env_theo}
\noindent
Now we can state our result for the uniform distribution in an ellipse:
\begin{env_cor}
\label{cor_ellipse_uniform}
Consider uniformly distributed points inside an ellipse $E$ with major axis $a=2$ and minor axis $2b<2$. If the  major axis is
placed between $(-1,0)$ and $(1,0)$, $E$ satisfies A1 - A6. Since the border functions of $E$ are given by $\pm b\sqrt{1-x^2}$, we get by symmetry $q_1 = \ldots = q_4 = 2b$.
Because of $\lambda^2(E) = \pi b$, we have $p_1 = \ldots = p_4 = 1/(\pi b)$. Hence, Theorem~\ref{satz_Vert_Konv_Hauptsatz} (and thus also Theorem~\ref{satz_Vert_Konv_M_n}) is applicable with i.i.d. random elements
$
\mathbf{Z}_{1},\ldots,\mathbf{Z}_{4} \sim  \NA_{\infty}\left( \sigma_1, \tau_1\right)
$.
\end{env_cor}

\begin{figure}[ht]
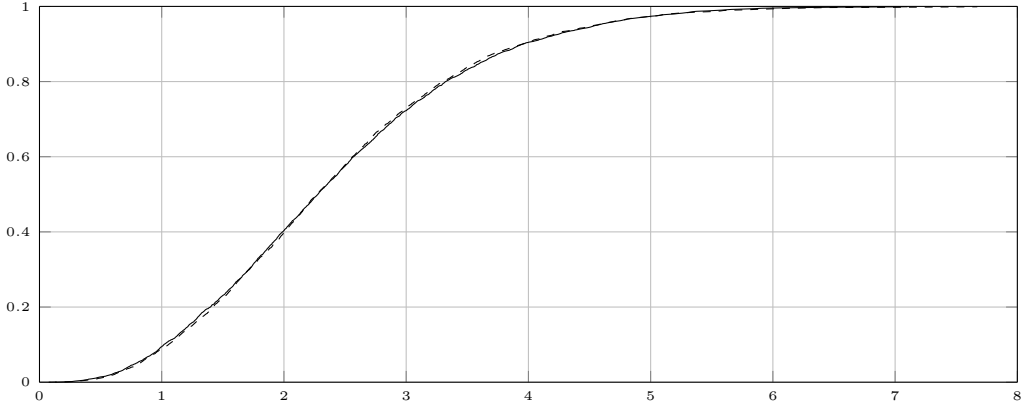

\begin{center}
%
\end{center}
	\caption{Empirical distribution function in the setting of Corollary~\ref{cor_ellipse_uniform} with $b=1/2$, $n=1000$ (solid, 5000 replications). The limit distribution is approximated as described in Corollary~\ref{cor_ellipse_uniform} for $m=8$ (dashed, 5000 replications).}
	\label{fig_cdf_plot}
\end{figure}

\cite{RaoJanson2015} described the limit distribution as that of
$$
\pi^{2/3} \min_{i,j \in \N}\left\{ x_i'+x_j'' - \frac{b^2}{4}(y_i'-y_j'')^2 \right\},
$$
where $\left\{ (x_i',y_i')\right\}$ and $\big{\{} (x_j'',y_j'')\big{\}}$ are two independent Poisson processes in the parabola $\left\{ (x,y) \in \R^2 :y^2 \le 2x\right\}$ with intensity 1. For simulations, the representation of the limit distribution given in \eqref{eq_Vert_Konv_Aussage} is much more useful, since the latter can easily be approximated. For the latter purpose, fix  $m \ge 1$ and replace $\min_{k,l \in \N}$ in the definition of $S^{i,j}$ in Theorem~\ref{satz_Vert_Konv_Hauptsatz} by $\min_{1 \le k,l \le m}$. This approximation is a consequence of Lemma~\ref{lemma_Nullfolge_Epsilon_k}. The bigger the minor half-axis $b$ is
(i.e. the more $E$ becomes 'circlelike'), the bigger $m$ has to be chosen in order to have a good approximation of the distributional limit in
 \eqref{eq_Vert_Konv_Aussage} (we omit details), see Figure~\ref{fig_cdf_plot}.

\section{Proof of Theorem~\ref{satz_Vert_Konv_Hauptsatz}}
\label{sec_proofs}

\subsection{The stochastic model}
\label{subsec_stoch_model}
Let $\widetilde N$ have a Poisson distribution Po$(n)$ and, independently of $N$, let  $Y_1,Y_2,\ldots $ be i.i.d. random variables with the same distribution as $X_1$.
Then $\Phi_n \overset{\mathcal{D}}{=} \widetilde{\Xi}_{\widetilde N} := \left\{ Y_1,Y_2,\ldots,Y_{\widetilde N}\right\}$ and hence
\begin{equation}
\label{eq_Vert_Gleich_Zwei_diam}
\diam(\Phi_n) \, \overset{\mathcal{D}}{=} \,  \diam(\widetilde{ \Xi}_{\widetilde N}).
\end{equation}
It will be useful to discriminate the points of $\widetilde{ \Xi}_{\widetilde N}$ according to the quadrant in which they realize.
To this end, put $ l_i := \int_{E_i}f(z)\mathrm{d}z$, $i \in \left\{ 1,2,3,4\right\}$.
For $n \ge 1$, let $N_i \sim \Po(n l_i)$ and, independently of $N_i$, let
$
X_1^{i},X_2^{i},\ldots $ be i.i.d. $\sim  l_i^{-1}\P_{X_1} \big|_{E_i},
$
where $\P_{X_1}\big|_{E_i}$ is the restriction of $\P_{X_1}$ to the set $E_i$.
The densities of these distributions are given by $f_i := l_i^{-1}f\big|_{E_i}$.
With $\Xi_{N_i}^{i} := \left\{ X_1^{i},X_2^{i},\ldots,X_{N_i}^{i}\right\}$, $N:= N_1 + \ldots + N_4$
and $ \Xi_N := \Xi_{N_1}^{1} \cup \Xi_{N_2}^{2} \cup \Xi_{N_3}^{3} \cup \Xi_{N_4}^{4}$
we get $\widetilde{ \Xi}_{\widetilde N} \overset{\mathcal{D}}{=}  \Xi_N$ and therefore
\begin{equation}
\label{eq_Vert_Gleich_zwei_Xi}
\diam(\widetilde{ \Xi}_{\widetilde N}) \, \overset{\mathcal{D}}{=} \,  \diam( \Xi_N).
\end{equation}
Because of this equality in distribution and \eqref{eq_Vert_Gleich_Zwei_diam}, it is sufficient to investigate $\diam( \Xi_N)$.
With the notation
$$
M_N^{j_1,j_2} := \max \left\{ \big| X_{k_1}^{j_1} - X_{k_2}^{j_2} \big|: 1 \le k_1 \le N_{j_1},1 \le k_2 \le N_{j_2}\right\}, \quad 1 \le j_1 \le j_2 \le 4,
$$
we have
$
\diam( \Xi_N) = \max_{1 \le j_1 \le j_2 \le 4} M_N^{j_1,j_2}.
$
The conditions A1 to A3 and A7 in mind, it is obvious that for sufficiently large $n$ only the pairs $ (1,2),(1,3),(2,4)$ and $(3,4)$ can be relevant for $(j_1,j_2)$. We obtain the important convergence
\begin{equation}
\label{eq_Stoch_Konv_diam_und_Max_von_4_max}
\P\bigg( \diam(\Xi_N) \neq \max\left\{ M_{N}^{1,2} ,M_{N}^{1,3} ,M_{N}^{2,4} ,M_{N}^{3,4}\right\}\bigg) \rightarrow  0.
\end{equation}
For example, $M_N^{1,2}$ will be determined by points inside the first (resp. second) quadrant, which are located close to the right (resp. left) pole of $E$ (cf. A1 to A3 and A7).
By A6, these points can easily be characterized by their distance to the origin: the norms of these points are close to the largest possible value $a$.
In Subsection~\ref{subsec_einzelner_Quad} we thus study  the asymptotic behaviour of those points with the largest norms. To this end, let
$
X_{(1)}^{i}, X_{(2)}^{i},\ldots, X_{(N_i)}^{i}
$
be a relabelling of  $X_1^{i},\ldots,X_{N_{i}}^{i}$ according to decreasing values of their norms, i.e., we have
$
\big|X_{(1)}^{i}\big| \ge \big|X_{(2)}^{i}\big| \ge \ldots \ge \big|X_{(N_i)}^{i}\big|.
$
The maximum interpoint distance of points with one of the $k$ largest norms inside one quadrant for fixed $k \ge 1$ will be denoted by
$$
M_{N,k}^{j_1,j_2} := \max \Big\{ \big| X_{(k_1)}^{j_1} - X_{(k_2)}^{j_2} \big|: 1 \le k_1 \le \min\left\{ k,N_{j_1}\right\},1 \le k_2 \le \min\left\{ k,N_{j_2}\right\}\Big\},
$$
for $(j_1,j_2) \in \big\{(1,2),(1,3),(2,4),(3,4)\big\}$.
Before we are able to characterize the behaviour of points lying close to a pole, we need some geometric results.


\subsection{Geometric considerations}
\label{subsec_Geom_Ueberl}
For $i \in \left\{ 1,2,3,4\right\}$ and small $h>0$ we set
\begin{equation}
\label{eq_Def_von_A_i_h}
  A_i(h) := E_i \backslash B(h).
\end{equation}
Because of A6, the interior of $A_i(h)$ is a connected set containing only points lying close to one of the poles.
It is easy to see that A7 implies
\begin{equation}
\label{eq_Verhalten_Integral_Ueber_Endkappe}
 \int_{A_i(h)}f(z)\mathrm dz \sim p_i\lambda^2\big( A_i(h)\big)
\end{equation}
for $i \in \left\{ 1,2,3,4\right\}$ as $h \to 0$, and for later usage we show the asymptotic behaviour of $\lambda^2\big( A_i(h)\big)$.
\begin{env_lemma}
\label{lemma_asymp_Verh_Endkappe_EllAehnlich}
For $i \in \left\{ 1,2,3,4\right\}$ and $c_i$ given in \eqref{def_der_Konstanten}, we have
$\lambda^2\big(A_i(h)\big)  \sim c_i  h^{3/2}$ as $h \to 0$.
\end{env_lemma}

\begin{proof}
W.l.o.g. let $i=1$, and fix $\varepsilon >0$.  From A5 we find a $x_0(\varepsilon) < a $ with
\begin{equation}
\label{eq_beweis_Flaeche_Ell_aehnlich_Anordnung}
(q_1-\varepsilon)f_a(x) \le g_1(x) \le (q_1 + \varepsilon)f_a(x)
\end{equation}
for each $x \in \big[ x_0(\varepsilon),a\big]$. Writing $B_h(x) := \sqrt{(a-h)^2-x^2}$ for the upper boundary function of $B(h)$
and $x^{\pm}(h)$ for the abscissae of the points of intersection between $(q_1\pm \varepsilon)f_a(x)$ and $B_h(x)$, some algebra gives
$$
x^{\pm}(h) = \sqrt{a^2-\frac{4 h (2a-h)}{4-(q_1\pm \varepsilon)^2}}.
$$
Putting $
  F^{\pm}(\varepsilon) := \left\{ (x,y) \in \R^2: 0 \le x \le a \text{ and } 0 \le y \le (q_1\pm\varepsilon)f_a(x)\right\}
$
and using \eqref{eq_beweis_Flaeche_Ell_aehnlich_Anordnung} we obtain
\begin{equation}
\label{eq_beweis_Flaeche_Ell_aehnlich_Anordnung_2}
\lambda^2\Big( F^{-}(\varepsilon) \backslash B(h)\Big)  \le \lambda^2\big( A_1(h)\big) \le \lambda^2\Big( F^{+}(\varepsilon) \backslash B(h)\Big).
\end{equation}
See figure~\ref{fig_Beweis_Flaeche} for an illustration.
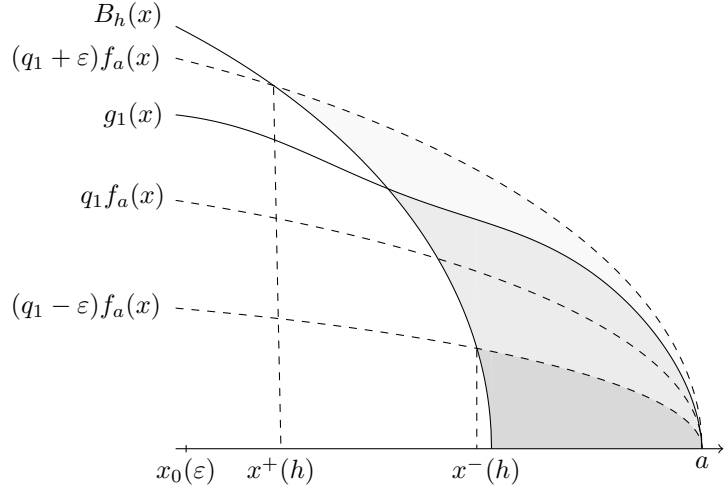
\begin{figure}[ht]
\begin{center}
\begin{tikzpicture}[xscale=14,yscale=12,domain=0.5:1.0,samples=300]
        \fill[lightgray!10] plot[domain=0.593:0.702] (\x,{0.5*sqrt(1-\x^2)})
            -- plot[domain=0.702:0.593]  (\x,{0.75*sqrt((1-0.2)^2-\x^2)});
        \fill[lightgray!10] plot[domain=0.702:1] (\x,{0.5*sqrt(1-\x^2)})
            -- plot[domain=1:0.70] (\x,{0.43*sqrt(1-\x^2)-0.01*cos(1000*\x) -0.01*sin(100*(\x+0.5))});
        \fill[lightgray!30] plot[domain=0.702:0.786] (\x,{0.43*sqrt(1-\x^2)-0.01*cos(1000*\x)
            -0.01*sin(100*(\x+0.5))}) -- plot[domain=0.786:0.702] (\x,{0.75*sqrt((1-0.2)^2-\x^2)});
        \fill[lightgray!30] plot[domain=0.786:1] (\x,{0.43*sqrt(1-\x^2)-0.01*cos(1000*\x)
            -0.01*sin(100*(\x+0.5))}) -- plot[domain=1:0.785] (\x,{0.18*sqrt(1-\x^2)});
        \fill[lightgray!60] plot[domain=0.786:0.8] (\x,{0.18*sqrt(1-\x^2}) -- plot[domain=0.8:0.786]
            (\x,{0.75*sqrt((1-0.2)^2-\x^2)});
        \fill[lightgray!60] plot[domain=0.8:1] (\x,{0.18*sqrt(1-\x^2}) -- plot[domain=1:0.799]
            (\x,{0});

        \draw[dashed] plot (\x,{0.318*sqrt(1-\x^2)})--(1,0);
        \draw[dashed] plot (\x,{0.5*sqrt(1-\x^2)})--(1,0);
        \draw[dashed] plot (\x,{0.18*sqrt(1-\x^2)})--(1,0);
        \draw[domain=0.5:1] plot (\x,{0.43*sqrt(1-\x^2)-0.01*cos(1000*\x) -0.01*sin(100*(\x+0.5))})--(1,0);

        \draw[domain = 0.5:0.8] plot (\x,{0.75*sqrt((1-0.2)^2-\x^2)})--(0.8,0);
        \draw[dashed] (0.593,0.401) -- (0.6,0) ;
        \draw[dashed] (0.786,0.11) -- (0.786,0);
        \draw (0.51,0.004) -- (0.51,-0.004);
        \draw (1,-0.015) node{$a$};
        \draw (0.51,-0.025) node{$x_0(\varepsilon)$};
        \draw (0.6,-0.025) node{$x^{+}(h)$};
        \draw (0.795,-0.025) node{$x^{-}(h)$};
        \draw (0.455,0.48) node{$B_h(x)$};
        \draw (0.415,0.433) node{$(q_1+\varepsilon)f_a(x)$};
        \draw (0.46,0.37) node{$g_1(x)$};
        \draw (0.45,0.28) node{$q_1f_a(x)$};
        \draw (0.415,0.158) node{$(q_1-\varepsilon)f_a(x)$};

       \draw[->] (0.5,0) -- (1.02,0);
\end{tikzpicture}
\end{center}
\caption{Illustration of the initial situation in the proof of Lemma~\ref{lemma_asymp_Verh_Endkappe_EllAehnlich} }
\label{fig_Beweis_Flaeche}
\end{figure}
To calculate the values on the left- and on the right-hand side of \eqref{eq_beweis_Flaeche_Ell_aehnlich_Anordnung_2} we use polar coordinates.
The upper boundary of $F^{\pm}(\varepsilon)$ is given by $\left\{ r^{\pm}(\varphi)(\cos \varphi,\sin \varphi): 0 \le \varphi \le \pi/2\right\}$ with
$$
 r^{\pm}(\varphi) := \frac{a (q_1\pm\varepsilon)}{2\sqrt{1+\left( \frac{(q_1\pm\varepsilon)^2 }{4}-1\right)\cos^2\varphi}}.
$$
If we denote by $\gamma^{\pm}_h$ the polar angle of the point $\Big( x^{\pm}(h),B_h\big( x^{\pm}(h)\big)\Big)$, it follows that
\begin{equation}
\label{eq_Bew_Verh_Flaeche_Tangens}
\tan(\gamma^{\pm}_h) = \frac{B_h\big( x^{\pm}(h)\big)}{ x^{\pm}(h)}.
\end{equation}
This leads to
\begin{align*}
  \lambda^2\Big(F^{\pm}(\varepsilon) \backslash B(h)\Big)
  &= \int\limits_{0}^{\gamma_h^{\pm}}\int\limits_{a-h}^{r^{\pm}(\varphi)}r \mathrm dr \mathrm d\varphi\\
  &= \frac{1}{2}\cdot \int\limits_{0}^{\gamma_h^{\pm}}\left( \frac{a^2 \cdot (q_1\pm\varepsilon)^2} {4+4\left( \frac{(q_1\pm\varepsilon)^2}{4}-1\right)\cos^2\varphi}-(a-h)^2\right)\mathrm d\varphi\\
  &= \frac{ a^2\cdot (q_1 \pm \varepsilon)}{4}\cdot \arctan \left( \frac{2}{q_1 \pm \varepsilon} \cdot \frac{B_h\big( x^{\pm}(h)\big)}{ x^{\pm}(h)}\right) - \frac{(a-h)^2}{2}\arctan\left(\frac{B_h\big( x^{\pm}(h)\big)}{ x^{\pm}(h)} \right).
\end{align*}
By expanding this term about $h=0$, we get
\begin{equation}
\label{eq_beweis_Flaeche_Ell_aehnlich_Entw_plus}
 \lambda^2\Big(F^{\pm}(\varepsilon) \backslash B(h)\Big) = c^{\pm}_1(\varepsilon)  h^{3/2} + O\left( h^{5/2}\right)
\end{equation}
with
$$
c^{\pm}_1(\varepsilon):= \frac{2(q_1 \pm \varepsilon)\sqrt{2a}}{3\sqrt{4-(q_1\pm\varepsilon)^2}}.
$$
Hence, by \eqref{eq_beweis_Flaeche_Ell_aehnlich_Anordnung_2} and \eqref{eq_beweis_Flaeche_Ell_aehnlich_Entw_plus} the inequality
$
\lambda^2\big( A_1(h)\big) \le c^{+}_1(\varepsilon)h^{3/2} + O(h^{5/2})
$
holds as $h \to 0$, and we have
  \begin{equation}
  \label{eq_beweis_Flaeche_vorletzer_Schritt}
  \frac{\lambda^2\big(A_1(h)\big) - c_1 h^{3/2}}{c_1 h^{3/2}}
  \le \frac{c^{+}_1(\varepsilon) h^{3/2} + O\left( h^{5/2}\right)-c_1 h^{3/2}}{c_1 h^{3/2}}
  = \frac{c^{+}_1(\varepsilon)-c_1}{c_1} + O\left( h\right) .
  \end{equation}
Now fix $\delta>0$. Since the function
$$
c^{+}_1(x) := \frac{2(q_1 + x)\sqrt{2a}}{3\sqrt{4-(q_1+x)^2}}
$$
is continuous from the right at $x=0$ for all valid $a$ and $q_1$, we can choose $\varepsilon >0$ in such a way that
$
(c^{+}_1(\varepsilon)-c_1)/c_1  \le \delta/2.
$
By \eqref{eq_beweis_Flaeche_vorletzer_Schritt} we get
$
\lambda^2\big(A_1(h)\big) - c_1  h^{3/2}  \le \delta c_1  h^{3/2}
$
as $h \to 0$. In the same way one can show
$
-\delta c_1 h^{3/2} \le \lambda^2\big(A_1(h)\big) - c_1 h^{3/2}
$ as $h \to 0$, and the proof is finished.
\end{proof}

Now, for sufficiently small $h$ (see A6), let $\eta_i(h)$ be the polar angle of the point
$
(\overline x_i(h), B_h\big( \overline x_i(h)\big)),
$
$i \in \left\{ 1,2\right\}$, and let $\eta_i(h)$ be the polar angle of
$
( \overline x_i(h), -B_h\big( \overline x_i(h)\big))
$
if $i \in \left\{ 3,4\right\}$. Furthermore, put
\begin{equation} \label{def_gamma_h_i}
\gamma_i(h) :=
\begin{cases}
  \eta_i(h), & i = 1,\\
  \pi - \eta_i(h), & i = 2,\\
  \eta_i(h) - \pi, & i = 3,\\
  2\pi - \eta_i(h), & i = 4.
\end{cases}
\end{equation}

\begin{env_lemma}
\label{lem_Verhalten_gamma_h_mit_Streckung_q_i}
For $i \in \left\{ 1,2,3,4\right\}$ we have $\gamma_i(h) \sim \tau_i  h^{1/2}$ as $h \to 0$, with $\tau_i$ given in \eqref{def_der_Konstanten}.
\end{env_lemma}
\begin{proof}
W.l.o.g. let $i=1$. Fix an arbitrary $\delta>0$. We have to show that the inequalities
$$
- \delta \le \frac{ \gamma_1(h) - \tau_1 h^{1/2}}{\tau_1  h^{1/2}} \le \delta
$$
hold for sufficiently small $h$. With the same notation as in the proof of Lemma~\ref{lemma_asymp_Verh_Endkappe_EllAehnlich} and \eqref{eq_Bew_Verh_Flaeche_Tangens}, it follows that
$$
\arctan\left( \frac{B_h\big( x^{-}(h)\big)}{ x^{-}(h)} \right) \le \gamma_1(h) \le \arctan\left( \frac{B_h\big( x^{+}(h)\big)}{ x^{+}(h)} \right).
$$
Power series expansions about $h=0$ show
$$
\arctan\left( \frac{B_h\big( x^{\pm}(h)\big)}{ x^{\pm}(h)} \right)
=\tau^{\pm}_1(\varepsilon)  h^{1/2} + O\left(h^{3/2}\right)
$$
with
$$
\tau^{\pm}_1(\varepsilon):=\frac{(q_1\pm \varepsilon)\sqrt{2a}}{a\sqrt{4-(q_1\pm\varepsilon)^2}}.
$$
The rest of the proof is by complete analogy with the proof of Lemma~\ref{lemma_asymp_Verh_Endkappe_EllAehnlich}.
\end{proof}

\noindent
For the sake of readability, we change our notation until the end of this subsection by using capitals for deterministic sequences. Moreover, we denote the underlying quadrant by a subscript instead of a superscript.
For $i \in \left\{ 1,2,3,4\right\}$ the value $N_{i,n}$ and $V_{i,n}$ denote the norm and the polar angle of the $n$-th deterministic point, respectively.
As in \eqref{def_gamma_h_i} we set
\begin{equation}
\label{eq_Umrechnung_V_i_n_zu_W_i_n}
W_{i,n} :=
\begin{cases}
  V_{i,n}, & i = 1,\\
  \pi - V_{i,n}, & i = 2,\\
  V_{i,n} - \pi, & i = 3,\\
  2\pi - V_{i,n}, & i = 4
\end{cases}
\end{equation}
and define a function $p:\R_{\ge 0}\times [0,2\pi) \to \R^2$ by $p(r,\varphi) :=  r(\cos \varphi, \sin \varphi)$.
For $i \in \left\{ 1,2,3,4\right\}$ we write
$
(N_{i,n},W_{i,n}) \rightarrow \Pol_i,
$
if $N_{i,n}\ge 0, W_{i,n} \ge 0$ for each $n$, the point $p(N_{i,n},V_{i,n})$ lies in $ E_i$ and for $i \in \left\{ 1,4\right\}$ (resp. $i \in \left\{ 2,3\right\}$) the points $p(N_{i,n},V_{i,n})$ converge to $ (a,0)$ (resp. $(-a,0)$).
Notice that $(N_{i,n},W_{i,n}) \rightarrow  \Pol_i$ implies $W_{i,n} \to 0$.

\begin{env_lemma}
\label{lemma_det_Dist_1_3}
Let $(N_{1,n},W_{1,n})$ and $(N_{3,n},W_{3,n})$ be deterministic sequences satisfying
$
(N_{i,n},W_{i,n}) \rightarrow   \Pol_i$, $i \in \{1,3\}.$
Then
\begin{align}
\;\big| p\left( N_{1,n},V_{1,n}\right) - p\left( N_{3,n},V_{3,n}\right) \big|
  &=\;N_{1,n} + N_{3,n} - \frac{a}{4}E_n^2 + \widetilde R_n, \label{eq_Resultat_TaylorEntw_Quadrant_1_3}
\end{align}
where
$$
  E_n := W_{1,n}-W_{3,n}, \ \
  \widetilde R_n  := O\left( E_n^4\right) + A_n + B_n + C_n + D_n
$$
and
\begin{eqnarray*}
  A_n & := & \frac{1}{4}\left(  \frac{1}{2}E_n^2 +O\left(E_n^4\right)\right)(a - N_{1,n}), \quad
  B_n := \frac{1}{4}\left(  \frac{1}{2}E_n^2 +O\left(E_n^4\right)\right)(a - N_{3,n}),\\
  C_n & := & -\frac{a}{16}\left( \frac{1}{2}E_n^2 +O\left(E_n^4\right)\right)^2, \quad
  D_n = O\left( (a-N_{1,n})^2+ (a-N_{3,n})^2 + \left( \frac{1}{2}E_n^2 + O\left( E_n^4\right)\right)^2\right).
\end{eqnarray*}
\end{env_lemma}
\begin{proof}
By the law of cosines and $V_{3,n} = \pi + W_{3,n}$ we get
\begin{align*}
  \big| p\left( N_{1,n},V_{1,n}\right) - p\left( N_{3,n},V_{3,n}\right) \big|^2
  &= \big| p\left( N_{1,n},W_{1,n}\right) - p\left( N_{3,n},\pi + W_{3,n}\right) \big|^2\\
  &= N_{1,n}^2 + N_{3,n}^2 - 2 N_{1,n}N_{3,n}\cos (\pi + W_{3,n} - W_{1,n})\\
  &= N_{1,n}^2 + N_{3,n}^2 + 2 N_{1,n}N_{3,n}\cos ( W_{1,n} - W_{3,n}).
\end{align*}
Using
$
\cos(x-y) = 1 - (x-y)^2/2 + O\big((x-y)^4\big)
$
as $(x,y) \to 0$ yields
\begin{align*}
  \big| p\left( N_{1,n},V_{1,n}\right) - p\left( N_{3,n},V_{3,n}\right) \big|
  =  \sqrt{ N_{1,n}^2 + N_{3,n}^2 + 2 N_{1,n}N_{3,n}\left( 1 - \frac{1}{2}E_n^2 + O\left( E_n^4\right)\right)}.
\end{align*}
Taylor's theorem for multivariate functions then gives
$$
 \sqrt{x^2+y^2+2xy\left( 1-z\right)} \ = \
x+y - \frac{a}{2}z + \frac{1}{4}z(a-x) + \frac{1}{4}z(a-y) - \frac{a}{16}z^2 +O\Big( (a-x)^2+(a-y)^2 +z^2\Big)
$$
as $(x,y,z) \rightarrow (a,a,0)$.
Putting $x = N_{1,n}$, $y = N_{3,n}$ and $z = \frac{1}{2}E_n^2 + O\left( E_n^4\right)$ leads to
$$
\big| p\left( N_{1,n},V_{1,n}\right) - p\left( N_{3,n},V_{3,n}\right) \big|
  = N_{1,n} + N_{3,n} - \frac{a}{2}\left( \frac{1}{2}E_n^2 +O\left( E_n^4\right)\right) + R_n \notag
  = N_{1,n} + N_{3,n} - \frac{a}{4}E_n^2 + \widetilde R_n.
$$
\end{proof}
\noindent
By the same reasoning, we have:
\begin{env_lemma}
\label{lemma_det_Dist_1_2}
Let $(N_{1,n},W_{1,n})$ and $(N_{2,n},W_{2,n})$ be deterministic sequences satisfying
$
(N_{i,n},W_{i,n}) \rightarrow   \Pol_i$, $i \in \{1,2\}.
$
Then
$
\big| p\left( N_{1,n},V_{1,n}\right) - p\left( N_{2,n},V_{2,n}\right) \big|
  = N_{1,n} + N_{2,n} - aF_n^2/4 + \widetilde R_n,
$
where
$
F_n:= W_{1,n}+W_{2,n}
$
instead of $E_n$ and $\widetilde R_n$ adjusted accordingly.
\end{env_lemma}
\noindent
Because of symmetry, Lemma~\ref{lemma_det_Dist_1_3} (resp. \ref{lemma_det_Dist_1_2}) can be applied to sequences in the second and fourth (resp. third and fourth) quadrant.


\subsection{A single quadrant}

\label{subsec_einzelner_Quad}
We now study the joint asymptotic behaviour of those points inside a fixed quadrant $Q_i$ that have the $k$ largest norms, where $k \ge 1 $ is fixed.
Since the  number $N_i$ of points follows a Poisson distribution, which means that $\P\left(N_i <k\right)>0$ for every $n \in \N$, we put
\begin{equation}
\label{def_Konvention_2}
X_{(j)}^{i} := 0
\textrm{ for } j \in \left\{ N_i+1,\ldots,k\right\} \textrm{ provided that } N_i < k.
\end{equation}
\noindent
We start with the following lemma, the proof of which is omitted.

\begin{env_lemma}
\label{lemma_Vert_Konv_Ord_Stat_mit_Poisson}
Let $N \sim \Po(\mu)$ and, independently of $N$, let $U_1,U_2, \ldots $ be i.i.d. $\sim \textrm{U}\big([0,1]\big)$. Writing
$
U_{(1)} \le U_{(2)} \le \ldots \le U_{(N)}
$
for the order statistics of $U_1,\ldots, U_N$, we have for fixed $k$
$$
\mu\left(U_{(1)}  , U_{(2)}  , \ldots  ,  U_{(k)}\right) \overset{\mathcal{D}}{\longrightarrow }\left( S_1  ,  S_2 ,  \ldots  , S_k\right)
$$
as $\mu \to \infty$. Here, $S_m = \sum_{j=1}^m Y_j$, and $Y_1,\ldots,Y_n$ are i.i.d. unit exponential random variables.
\end{env_lemma}

\begin{env_lemma}
\label{lem_Vert_Konv_Normen}
Let $k \ge 1$ be fixed. Based on $Y_1,\ldots,Y_k$ and $S_m$ as above, we have
$$
n^{2/3}
\big(
a- \big|X_{(1)}^{i}\big|,  \ldots  ,  a- \big|X_{(k)}^{i}\big|
\big)
\overset{\mathcal{D}}{\longrightarrow }
\sigma_i\left( S_1^{2/3}  , S_2^{2/3}  ,  \ldots  , S_k^{2/3}\right),
$$
where $\sigma_i $ is given in \eqref{def_der_Konstanten}.
\end{env_lemma}
\begin{proof}
Since $X_1^{i}$ has a $\lambda^2$-density, the distribution function $F$ of $\left|X_1^{i}\right|$ is continuous and hence $1-F\big(\left|X_1^{i}\right|\big) \sim \textrm{U}\big([0,1]\big).$ Independently of $N_i$ and $X_{(1)}^{i},\ldots,X_{(N_i)}^{i}$, let $U_1,U_2, \ldots $ be i.i.d. $\sim \textrm{U}\big([0,1]\big)$ and write
$
U_{(1)} \le U_{(2)} \le \ldots \le U_{(N_i)}
$
for the order statistics of $U_1,\ldots, U_{N_i}$.
Because of $\P\left( N_i < k\right) \rightarrow  0$, the case $\left\{ N_i < k\right\}$ is negligible in the following.
We get
$$
\big(
1-F\big(\big|X_{(1)}^{i}\big|\big) ,  \ldots  ,  1-F\big(\big|X_{(k)}^{i}\big|\big)
\big) \sim
\big( U_{(1)}  ,  U_{(2)}  ,  \ldots  ,  U_{(k)}\big)
$$
and, due to $N_i \sim \Po(n l_i)$, Lemma~\ref{lemma_Vert_Konv_Ord_Stat_mit_Poisson} yields
$
n  l_i\big(
1-F\big(\big|X_{(1)}^{i}\big|\big) ,  \ldots  ,  1-F\big(\big|X_{(k)}^{i}\big|\big)\big)
\overset{\mathcal{D}}{\longrightarrow } \big(S_1  ,  \ldots  ,  S_k\big).
$
Assumption A6, the equations~\eqref{eq_Def_von_A_i_h} and \eqref{eq_Verhalten_Integral_Ueber_Endkappe} and Lemma~\ref{lemma_asymp_Verh_Endkappe_EllAehnlich} show
\begin{align*}
  1 - F(x)
  &= \P\big(\left|X_1^{i}\right| > x \big)
  = \int_{A_i(a-x)}\frac{f(z)}{ l_i}\mathrm{d}z
  \sim \frac{p_i}{ l_i} \lambda^2\big( A_i(a-x)\big)
  \sim \frac{p_i}{ l_i}c_i(a-x)^{3/2}
\end{align*}
as $x \to a$.
Therefore,
$
nl_i\big(
1-F\big(\big|X_{(1)}^{i}\big|\big), \ldots , 1-F\big(\big|X_{(k)}^{i}\big|\big)
\big)
$
and
$
n p_i c_i \big(
\big(a-\big|X_{(1)}^{i}\big|\big)^{3/2} ,  \ldots ,  \big(a-\big|X_{(k)}^{i}\big|\big)^{3/2}
\big)
$
are asymptotically equivalent and thus have the same limit $\big(S_1 ,  \ldots  , S_k \big)$ in distribution. With the function
$
g(x_1,\ldots,x_k):= \big( x_1^{2/3},\ldots,x_k^{2/3}\big)
$
and the continuous mapping theorem we obtain
\begin{eqnarray*}
  n^{2/3}
\big(a- \big|X_{(1)}^{i}\big|,  \ldots  ,  a- \big|X_{(k)}^{i}\big|\big) & = &
\left( \frac{1}{p_ic_i} \right)^{2/3} \left( n p_i c_i\right)^{2/3}
\big(a- \big|X_{(1)}^{i}\big|,  \ldots  ,  a- \big|X_{(k)}^{i}\big|\big)\\
& = & \sigma_i g\left( n p_i c_i \big(a-\big|X_{(1)}^{i}\big|\big)^{3/2}  ,  \ldots  , n p_i c_i\big(a-\big|X_{(k)}^{i}\big|\big)^{3/2} \right)\\
& \overset{\mathcal{D}}{\rightarrow } &  \sigma_ig \left( S_1  ,  S_2  ,  \ldots  ,  S_k\right)\\
& = &  \sigma_i \big( S_1^{2/3}  ,  S_2^{2/3}  ,  \ldots  ,  S_k^{2/3}\big).
\end{eqnarray*}
\end{proof}
In what follows, observe that $k$ is fixed in this subsection.
If $k \le N_i$, let $\eta_{(l)}^{i}$ the polar angle of $X_{(l)}^{i}$ for $l=1,\ldots,k$, and set $\eta_{(l)}^{i} := 0$ for $l \in \left\{ N_i+1,\ldots,k\right\}$ if $ N_i < k$.
Furthermore, put
\begin{equation}
\label{def_Polarwinkel_zuf_Punkte}
\gamma_{(l)}^{i} :=
\begin{cases}
  \eta_{(l)}^{i}, & i = 1,\\
  \pi - \eta_{(l)}^{i}, & i = 2,\\
  \eta_{(l)}^{i} - \pi, & i = 3,\\
  2\pi - \eta_{(l)}^{i}, & i = 4,
\end{cases},
\qquad l \in \left\{ 1,\ldots,k\right\}.
\end{equation}

We need the joint asymptotic behaviour of the angles $\gamma_{(1)}^{i},\ldots,\gamma_{(k)}^{i}$. As before, the case $N_i <k$ is negligible.
In Lemma~\ref{lem_Vert_Konv_Normen} we have shown the weak convergence of $n^{2/3}\big( a - \big|X_{(j)}^{i}\big| \big)$ for $j \in \left\{ 1,\ldots,k\right\}$. Let $0 < x_1 < \ldots < x_k$
 and consider  the conditions
$
n^{2/3}\big(a-  \big|X_{(j)}^{i}\big|\big) = x_j >0,
$
$j=1,\ldots, k$.
Under these conditions, the angles $\gamma_{(1)}^{i},\ldots,\gamma_{(k)}^{i}$ are asymptotically independent.
Some calculations and Lemma~\ref{lem_Verhalten_gamma_h_mit_Streckung_q_i} show that the conditional density of $n^{1/3}\gamma_{(j)}^{i}$ converges pointwise to the density of a U$\left( \big[ 0,\tau_i\sqrt{x_j}\right]\big)$-distributed random variable
(with $\tau_i$ given in \eqref{def_der_Konstanten}).
We thus have the following result:

\begin{env_lemma}
\label{lemma_Gem_Konv_Norm_Winkel_k}
      For fixed $k \ge 1$ let $Y_1,\ldots,Y_k$ and $U_1,\ldots,U_k$ be independent random variables, where $Y_1,\ldots,Y_k$ are i.i.d. with a unit exponential distribution, and $U_1,\ldots,U_k$ are i.i.d. with the uniform distribution U$([0,1])$. For $i\in \{1,2,3,4\}$, let
$Z_{1,m}^i \, :=  \, \sigma_i\left(\sum_{j=1}^{m}Y_j\right)^{2/3}$ and $Z_{2,m}^i \, := \, U_m \tau_i(Z_{1,m}^i)^{1/2}$ for $m=1,\ldots,k$. Then
$$
\left(
 n^{2/3}\big( a-  \big|X_{(1)}^{i}\big| \big) ,   n^{1/3}\gamma_{(1)}^{i}  , \ldots ,
    n^{2/3}\big( a-  \big|X_{(k)}^{i}\big| \big) ,   n^{1/3}\gamma_{(k)}^{i}
  \right) \overset{\mathcal{D}}{\longrightarrow } \left(
      Z_{1,1}^i , Z_{2,1}^i ,\ldots , Z_{1,k}^i, Z_{2,k}^i
  \right).
$$
\end{env_lemma}
The limit distribution in Lemma~\ref{lemma_Gem_Konv_Norm_Winkel_k} shall be denoted by $\NA_k(\sigma_i,\tau_i)$ ('norm-angle distribution of order $k$').


\subsection{Different quadrants}
\label{subsec_mehrere_Quad}
Instead of studying $\diam({\Xi}_N)$, we make two restrictions at this point: On the one hand we examine the behaviour of the diameters $M_N^{j_1,j_2}$ separately for each pair $(j_1,j_2) \in \left\{ (1,2),(1,3),(2,4),(3,4)\right\}$. On the other hand, we approximate these random variables by $M_{N,k}^{j_1,j_2}$. In a first step we fix $k \ge 1$, and then we let $k$ tend to infinity. We will see, that the difference between $M_N^{j_1,j_2}$ and $M_{N,k}^{j_1,j_2}$ becomes asymptotically negligible.
For fixed $k$ and $i \in \left\{ 1,2,3,4\right\}$ we have
$
\big( \big| X_{(l)}^{i}\big|, \gamma_{(l)}^{i}\big) \rightarrow \Pol_i
$
almost surely
for each $l \in \left\{ 1,\ldots,k\right\}$. In other words, the series expansions of Lemmas~\ref{lemma_det_Dist_1_3} and \ref{lemma_det_Dist_1_2} can be applied to each pair
 $\big( \big| X_{(l)}^{i}\big|, \gamma_{(l)}^{i}\big)$.

\begin{env_prop}
\label{satz_Vert_Konv_1_3_je_k_groesste_Normen}
For fixed $k \ge 1$ we have
$
n^{2/3}\big( 2a - M_{N,k}^{1,3}\big)
\overset{\mathcal{D}}{\longrightarrow } \min_{1\le i,j \le k}\big{\{}Z_{1,i}^{1} + Z_{1,j}^{3} +\frac{a}{4}\left(Z_{2,i}^{1} - Z_{2,j}^{3}\right)^2 \big{\}}
$
with two independent random elements
$(Z_{1,j}^{1},Z_{2,j}^{1})_{j \ge 1} \sim \NA_{k}(\sigma_1,\tau_1)$ and
$(Z_{1,j}^{3},Z_{2,j}^{3})_{j \ge 1} \sim \NA_{k}(\sigma_3,\tau_3)$.
\end{env_prop}
\begin{proof}
As $n$ tends to infinity, the probabilities $\P\left( N_1 < k\right)$ and $\P\left( N_3 < k\right)$ converge to $0$. Because of this asymptotic negligibility, we assume $k < N_1$ and $ k < N_3$.
Since the points in $Q_1$ and $Q_3$ are independent, Lemma~\ref{lemma_Gem_Konv_Norm_Winkel_k} implies
\begin{equation}
\label{eq_bew_Vert_Konv_M_N_k_VertKonv_in_Quad}
\left(
  \begin{array}{ccccc}
    n^{2/3}\big( a-  \big|X_{(1)}^{1}\big| \big) , &  n^{1/3}\gamma_{(1)}^{1} & , \ldots , &
    n^{2/3}\big( a-  \big|X_{(k)}^{1}\big| \big) , &  n^{1/3}\gamma_{(k)}^{1} \\
    n^{2/3}\big( a-  \big|X_{(1)}^{3}\big| \big) , &  n^{1/3}\gamma_{(1)}^{3} & , \ldots , &
    n^{2/3}\big( a-  \big|X_{(k)}^{3}\big| \big) , &  n^{1/3}\gamma_{(k)}^{3}
  \end{array}
\right) \overset{\mathcal{D}}{\longrightarrow }
\left(
  \begin{array}{c}
    Z_{1}\\
    Z_{3}\\
  \end{array}
\right),
\end{equation}
with two independent random elements $Z_{1} \sim \NA_{k}(\sigma_1,\tau_1) $ and $Z_{3}\sim  \NA_{k}(\sigma_3,\tau_3)$.
Since
\begin{equation}
n^{2/3}\big( 2a - M_{N,k}^{1,3}\big)
= n^{2/3}\big( 2a - \max_{1 \le i,j \le k} \big|X_{(i)}^{1}-X_{(j)}^{3} \big|\big)
=  \min_{1 \le i,j \le k}\Big\{ n^{2/3}\big(2a -  \big|X_{(i)}^{1}-X_{(j)}^{3} \big| \big) \Big\}, \label{eq_bew_Vert_Konv_M_N_k_Umgeformt}
\end{equation}
we define
$
h^{-}\left(
 y_1  , y_2 , z_1 , z_2
\right) := y_1 + z_1 + a(y_2 - z_2)^2/4
$
and, using Lemma~\ref{lemma_det_Dist_1_3} for $(i,j) \in \left\{ 1,\ldots,k\right\}^2$, obtain
\begin{eqnarray*}
   n^{2/3}\big( 2a - \big|X_{(i)}^{1}-X_{(j)}^{3}\big|\big) & = &
   n^{2/3} \big( 2a - \big( \big|X_{(i)}^{1}\big| + \big|X_{(j)}^{3}\big| - \frac{a}{4}\left( \gamma_{(i)}^{1}-\gamma_{(j)}^{3}\right)^2 + \widetilde R_n  \big) \big)\\
& = &  n^{2/3} \big( a-  \big|X_{(i)}^{1}\big| \big) + n^{2/3} \big( a-  \big|X_{(j)}^{3}\big| \big) + \frac{a}{4} \left( n^{1/3}\gamma_{(i)}^{1}-n^{1/3}\gamma_{(j)}^{3}\right)^2 + n^{2/3}\widetilde{R}_n\\
 & =&  h^{-}\left(
    \begin{array}{cccc}
      n^{2/3} \Big( a-  \big|X_{(i)}^{1}\big| \Big) \, , \,
      n^{1/3}\gamma_{(i)}^{1} \, , \,
      n^{2/3} \Big( a-  \big|X_{(j)}^{3}\big| \Big) \, , \,
      n^{1/3}\gamma_{(j)}^{3}
    \end{array}
  \right)+n^{2/3}\widetilde R_n.
\end{eqnarray*}
To show that $n^{2/3}\widetilde R_n = o_{\P}(1)$ put $ E_n := \gamma_{(i)}^{1} -\gamma_{(j)}^{3}$. From  \eqref{eq_Resultat_TaylorEntw_Quadrant_1_3} we have
$
\widetilde R_n = O\left(E_n^4\right) + A_n + B_n + C_n + D_n,
$
where
 \begin{eqnarray*}
  A_n  & := & \frac{1}{4}\left( \frac{1}{2}E_n^2 + O\left( E_n^4\right)\right) \cdot
  \big( a - \big|X_{(i)}^{1}\big|\big), \quad
  B_n := \frac{1}{4}\left( \frac{1}{2}E_n^2 + O\left( E_n^4\right)\right) \cdot
  \big( a - \big|X_{(j)}^{3}\big|\big),\\
  C_n & := & -\frac{a}{16}  \left( \frac{1}{2}E_n^2 + O\left( E_n^4\right)\right)^2, \quad
    D_n  =  O\left( \big( a-\big|X_{(i)}^{1}\big|\big)^2 + \big( a - \big|X_{(j)}^{3}\big|\big)^2 + \bigg( \frac{1}{2}E_n^2 + O\left(E_n^4\right)\bigg)^2\right).
\end{eqnarray*}
Since $E_n = o_{\P}(1)$ and $(n^{1/3}E_n)$ is a tight sequence, we get $n^{2/3}O\left( E_n^4\right) = \left( n^{1/3}E_n\right)^2 O\left( E_n^2\right)= o_{\P}(1)$.
From
$
n^{2/3}A_n = \frac{1}{4}\big( E_n^2/2 + O\left( E_n^4\right)\big) \cdot \big (n^{2/3}
  \big( a - \big|X_{(i)}^{1}\big|\big) \big),
$
Lemma~\ref{lem_Vert_Konv_Normen} and $E_n = o_{\P}(1)$ we obtain $n^{2/3}A_n =o_{\P}(1)$. The same reasoning gives $n^{2/3}B_n = o_{\P}(1)$. Now,
$$
n^{2/3}C_n = -\frac{a}{16}  \bigg( \frac{1}{2}n^{1/3}E_n^2 + n^{1/3}O\left( E_n^4\right)\bigg)^2
= -\frac{a}{16}  \bigg( \frac{1}{2}\left( n^{1/3}E_n \right) E_n + \left( n^{1/3}E_n \right) O\left( E_n^3\right)\bigg)^2
$$
entails $n^{2/3}C_n = o_{\P}(1)$, since $E_n = o_{\P}(1)$ and $(n^{1/3}E_n)$ is tight. Lemma~\ref{lem_Vert_Konv_Normen} yields
$
n^{2/3}\big( a-\big|X_{(i)}^{m}\big|\big)^2 = o_{\P}(1)
$
for $m \in \left\{ 1,3\right\}$ and, together with $n^{2/3}C_n = o_{\P}(1)$, we get $n^{2/3}D_n = o_{\P}(1)$, i.e $n^{2/3}\widetilde R_n = o_{\P}(1)$.
We thus can rewrite \eqref{eq_bew_Vert_Konv_M_N_k_Umgeformt} as
\begin{align*}
n^{2/3}\big( 2a - M_{N,k}^{1,3}\big)
&=  \min_{1 \le i,j \le k}\left\{  h^{-}\left(
    \begin{array}{c}
      n^{2/3} \Big( a-  \big|X_{(i)}^{1}\big| \Big), \\
      n^{1/3}\gamma_{(i)}^{1}, \\
      n^{2/3} \Big( a-  \big|X_{(j)}^{3}\big| \Big), \\
      n^{1/3}\gamma_{(j)}^{3} \\
    \end{array}
  \right)
+ o_{\P}(1)\right\},
\end{align*}
and the assertion follows from the continuous mapping theorem and \eqref{eq_bew_Vert_Konv_M_N_k_VertKonv_in_Quad}.
\end{proof}

\noindent
The next lemma shows that the difference between $M_{N}^{1,3}$ and $M_{N,k}^{1,3}$ becomes negligible as $k$ tends to infinity.

\begin{env_lemma}
\label{lemma_Nullfolge_Epsilon_k}
There exists a null sequence $(\varepsilon_k)_{k\ge 1}$ with
$
\P\big(M_{N,k}^{1,3} \ne M_N^{1,3}\big) \le \varepsilon_k
$
for every $k\ge 1$ and sufficiently large $n$.
\end{env_lemma}
\begin{proof}Let  $\mathbf{B}^1_{N_1}$  be the $\R^\N$-valued random element with components
$n^{2/3} \big( a-  |X_{(1)}^{1}| \big) $,
$n^{1/3}\gamma_{(1)}^{1}$,$\ldots$,
$n^{2/3} \big( a-  |X_{(N_1)}^{1}| \big)$,
$n^{1/3}\gamma_{(N_1)}^{1}$,
followed by
$n^{2/3} \big( a-  \big|X_{(N_1)}^{1}\big| \big)$ and
$n^{1/3}\gamma_{(N_1)}^{1}$, repeated infinitely often.
Let $Y_1^{1},Y_2^{1},\ldots$ and $U_1^{1},U_2^{1},\ldots$ be independent random variables, where $Y_1^{1},Y_2^{1},\ldots$ are i.i.d. with a unit exponential distribution, and $U_1^{1},U_2^{1},\ldots$ are i.i.d. with the uniform distribution U$([0,1])$. For every $m \in \N$ we put $S_m^{1} :=  \sigma_1\left(\sum_{j=1}^{m}Y_j^{1}\right)^{2/3}$, $T_m^{1} :=  U_m^{1} \tau_1\sqrt{S_m^{1}}$ and finally
we set
$
\mathbf{R}_{1} := \left(S_j^{1},T_j^{1}\right)_{j\ge 1}.
$
From Lemma~\ref{lemma_Gem_Konv_Norm_Winkel_k} we obtain
$
\pi_k \left( \mathbf{B}_{N_1}^{1} \right)  \overset{\mathcal{D}}{\longrightarrow } \pi_k\left( \mathbf{R}_{1} \right),
$
for each fixed $k$, where $\pi_k: \R^\N \rightarrow \R^k$
denotes the projection onto the first $k$ components.
Since the class of finite-dimensional sets is a convergence-determining class for $\R^{\N}$ (see Example 2.4 in \cite{Billingsley1999}), we get
\begin{equation}
\label{eq_R_unendlich_Konv_erster_Quadrant}
\mathbf{B}_{N_1}^{1} \overset{\mathcal{D}}{\longrightarrow } \mathbf{R}_{1}.
\end{equation}
With similar definitions, we also conclude that
\begin{equation}
\label{eq_R_unendlich_Konv_dritter_Quadrant}
\mathbf{B}_{N_3}^{3} \overset{\mathcal{D}}{\longrightarrow } \mathbf{R}_{3}.
\end{equation}
Since the points in the first and in the third quadrant are independent, the limit distributions $\mathbf{R}_1$ and $\mathbf{R}_3$ are also independent. We now assume $k \le \min\left\{ N_1,N_3\right\}$, and
for $k \ge 1$ fixed and $1 \le i,j \le k$ we define
$$
p_{i,j,n} :=
  \P\left( \big|X_{(i)}^{1} - X_{(j)}^{3}\big| = \max\limits_{\substack{1 \le l \le N_1\\ 1 \le m \le  N_3}}\big|X_{(l)}^{1} - X_{(m)}^{3}\big|\right).
$$
With
$$
h_{i,j}^-: \begin{cases}
\R^\N \times \R^\N \rightarrow \R, \\ (x,y) \mapsto x_{2i-1} + y_{2j-1} + \frac{a}{4}\left( x_{2i}- y_{2j}\right)^2
\end{cases}
$$
and the same reasoning as in the proof of Proposition~\ref{satz_Vert_Konv_1_3_je_k_groesste_Normen} we get
\begin{align*}
  p_{i,j,n} &= \P\left( \big|X_{(i)}^{1} - X_{(j)}^{3}\big| =  \max\limits_{\substack{1 \le l \le N_1\\ 1 \le m \le  N_3}}\big|X_{(l)}^{1} - X_{(m)}^{3}\big|\right)\\
  &= \P\left( n^{2/3} \bigg(2a - \big|X_{(i)}^{1} - X_{(j)}^{3}\big|\bigg) = n^{2/3} \left(2a- \max\limits_{\substack{1 \le l \le N_1\\ 1 \le m \le  N_3}}\big|X_{(l)}^{1} - X_{(m)}^{3}\big|\right)\right)\\
  &= \P\left( n^{2/3} \bigg(2a - \big|X_{(i)}^{1} - X_{(j)}^{3}\big|\bigg) = \min\limits_{\substack{1 \le l \le N_1\\ 1 \le m \le  N_3}} \left\{n^{2/3} \bigg(2a-\big|X_{(l)}^{1} - X_{(m)}^{3}\big|\bigg) \right\}\right) \\
  &= \P \left ( h_{i,j}^{-}\left( \mathbf{B}_{N_1}^{1},\mathbf{B}_{N_3}^{3}\right) + o_{i,j} = \min\limits_{\substack{1 \le l \le N_1\\ 1 \le m \le  N_3}} \left\{ h_{l,m}^{-}\left( \mathbf{B}_{N_1}^{1},\mathbf{B}_{N_3}^{3}\right) + o_{l,m} \right\}\right),
\end{align*}
with $o_{i,j} = o_{\P}(1)$ and $o_{l,m} = o_{\P}(1)$. These stochastic sequences are written explicitly, since it will be important  that $o_{i,j} = o_{l,m}$ for $(i,j) = (l,m)$.
In view of the definition of $\mathbf{B}_{N_1}^{1}$ and $\mathbf{B}_{N_3}^{3}$ it is obvious that we can take the minimum over $l,m \in \N$ instead of $1 \le l \le N_1, 1 \le m \le  N_3$. Equations~\eqref{eq_R_unendlich_Konv_erster_Quadrant} and \eqref{eq_R_unendlich_Konv_dritter_Quadrant} and the continuous mapping theorem yield
$$
  p_{i,j,n}  = \P \bigg ( h_{i,j}^{-}\left( \mathbf{B}_{N_1}^{1},\mathbf{B}_{N_3}^{3}\right) + o_{i,j} = \min\limits_{l,m \in \N} \left\{ h_{l,m}^{-}\left( \mathbf{B}_{N_1}^{1},\mathbf{B}_{N_3}^{3}\right) + o_{l,m} \right\}\bigg)\rightarrow  p_{i,j},
$$
with
$$
p_{i,j} := \P \bigg ( h_{i,j}^{-}\left( \mathbf{R}_1,\mathbf{R}_3\right)  = \min\limits_{l,m \in \N} h_{l,m}^{-}\left( \mathbf{R}_1,\mathbf{R}_3\right) \bigg).
$$
Putting
$
B_n := \big\{ k \le \min\left\{ N_1,N_3\right\}\big\}
$
for $k \ge 1$ fixed, we get $\P(B_n) \rightarrow  1$ and consequently
$$
\begin{aligned}
  &\P\left(M_{N,k}^{1,3} \ne M_N^{1,3}\right)&& = && \P\left(M_{N,k}^{1,3} \ne M_N^{1,3}\big| B_n\right)\cdot \P(B_n) + \P\left(M_{N,k}^{1,3} \ne M_N^{1,3}\big| B_n^c\right)\cdot\P(B_n^c) \\
  & && = && \P\left(M_{N,k}^{1,3} \ne M_N^{1,3}\big| B_n\right)\big( 1 +o(1)\big) + o(1) \\
  & && = && \bigg(1 - \P\left( M_{N,k}^{1,3} = M_N^{1,3}\big| B_n\right)\bigg)\big( 1 +o(1)\big) +o(1)\\
  & && = && \bigg(1 -  \sum_{i,j=1}^{k}p_{i,j,n}\bigg)\big( 1 +o(1)\big) +o(1) \\
  & &&\hspace{-2mm}\overset{n\rightarrow \infty}{\longrightarrow } &&1 - \sum_{i,j=1}^{k}p_{i,j}.
\end{aligned}
$$
Since $\sum_{i,j=1}^{\infty}p_{i,j} = 1$ and $p_{i,j} \ge 0$, the probability above converges to $0$ as $k \rightarrow \infty$.
\end{proof}

\begin{env_prop}
\label{satz_Vert_Konv_1_3}
We have
\begin{equation}\label{konvergenz13}
n^{2/3}\big( 2a - M_{N}^{1,3}\big)
\overset{\mathcal{D}}{\longrightarrow } \min_{i,j \in \N}\left\{Z_{1,i}^{1} + Z_{1,j}^{3} +\frac{a}{4}\left(Z_{2,i}^{1} - Z_{2,j}^{3}\right)^2 \right\}
\end{equation}
with two independent random elements
$(Z_{1,j}^{1},Z_{2,j}^{1})_{j \ge 1}  \sim \NA_{\infty}(\sigma_1,\tau_1)$ and $(Z_{1,j}^{3},Z_{2,j}^{3})_{j \ge 1}\sim \NA_{\infty}(\sigma_3,\tau_3).$
\end{env_prop}
\begin{proof}
Write $F_n$ for the distribution function (df) of $n^{2/3}(2a- M_N^{1,3})$, and let $G$ be the df of the limit occurring in \eqref{konvergenz13}.
Furthermore, $G_k$ denotes the df of the right-hand side of \eqref{konvergenz13} with $\min_{k,l\in \N}$ replaced by $\min_{1\le i,j \le k}$.
For $k \rightarrow \infty$ we have
$$
\min\limits_{ 1 \le i,j \le k}\left\{Z_{1,i}^{1} + Z_{1,j}^{3} +\frac{a}{4}\left(Z_{2,i}^{1} - Z_{2,j}^{3}\right)^2 \right\}  \overset{\P}{\longrightarrow }
\min\limits_{ i,j \in \N}\left\{Z_{1,i}^{1} + Z_{1,j}^{3} +\frac{a}{4}\left(Z_{2,i}^{1} - Z_{2,j}^{3}\right)^2 \right\},
$$
and hence $G_k \overset{\mathcal{D}}{\longrightarrow } G$. Fix $t > 0$.
On the one hand, Proposition~\ref{satz_Vert_Konv_1_3_je_k_groesste_Normen} and Lemma~\ref{lemma_Nullfolge_Epsilon_k} shows
\begin{eqnarray*}
F_{n}(t) & = & \P\big( n^{2/3}\big( 2a - M_{N}^{1,3}\big) \le t\big)\\
 & = &\P\big( n^{2/3}\big( 2a - M_{N}^{1,3}\big)  \le t , M_{N,k}^{1,3} = M_N^{1,3}\big) +  \P\big( n^{2/3}\big( 2a - M_{N}^{1,3}\big)  \le t , M_{N,k}^{1,3} \ne M_N^{1,3}\big)\\
 &\le  &\P\big( n^{2/3}\big( 2a - M_{N}^{1,3}\big)  \le t \big| M_{N,k}^{1,3} = M_N^{1,3}\big) \cdot \P\left(M_{N,k}^{1,3} = M_N^{1,3}\right) + \varepsilon_k\\
 &\le & \P\big( n^{2/3}\big( 2a - M_{N}^{1,3}\big)  \le t \big| M_{N,k}^{1,3} = M_N^{1,3}\bigg)  + \varepsilon_k\\
 & =  &\P\big( n^{2/3}\big( 2a - M_{N,k}^{1,3}\big)  \le t \big) + \varepsilon_k\\
 & \overset{n\rightarrow \infty}{\longrightarrow }&  G_k(t) + \varepsilon_k.
\end{eqnarray*}
On the other hand we get
\begin{eqnarray*}
F_{n}(t) &=  &\P\big( n^{2/3}\big( 2a - M_{N}^{1,3}\big)  \le t , M_{N,k}^{1,3} = M_N^{1,3}\big) +     \P\big( n^{2/3}\big( 2a - M_{N}^{1,3}\big)  \le t , M_{N,k}^{1,3} \ne M_N^{1,3}\big)\\
 &\ge & \P\big( n^{2/3}\big( 2a - M_{N}^{1,3}\big)  \le t \big| M_{N,k}^{1,3} = M_N^{1,3}\big) \cdot \P\left(M_{N,k}^{1,3} = M_N^{1,3}\right)\phantom{ + \varepsilon_k}\\
&\ge &\P\big( n^{2/3}\big( 2a - M_{N,k}^{1,3}\big) \le t \big) \cdot(1-\varepsilon_k)\\
  & \overset{n\rightarrow \infty}{\longrightarrow } & G_k(t)\cdot(1-\varepsilon_k).
\end{eqnarray*}
Since $G_k(t) \rightarrow G(t)$ and $\varepsilon_k \rightarrow  0$ as $k \to \infty$ (see Lemma~\ref{lemma_Nullfolge_Epsilon_k}), the assertion follows.
\end{proof}
\noindent
For reasons of symmetry we get:
\begin{env_prop}
\label{satz_Vert_Konv_2_4}
We have $
n^{2/3}\big( 2a - M_{N}^{2,4}\big)
\overset{\mathcal{D}}{\longrightarrow } \min_{i,j \in \N}\left\{Z_{1,i}^{2} + Z_{1,j}^{4} +\frac{a}{4}\left(Z_{2,i}^{2} - Z_{2,j}^{4}\right)^2 \right\}
$
with two independent random elements
$(Z_{1,j}^{2},Z_{2,j}^{2})_{j \ge 1}  \sim \NA_{\infty}(\sigma_2,\tau_2)$ and $(Z_{1,j}^{4},Z_{2,j}^{4})_{j \ge 1}\sim \NA_{\infty}(\sigma_4,\tau_4).$
\end{env_prop}
\noindent
The same steps as above and Lemma~\ref{lemma_det_Dist_1_2} yield:
\begin{env_prop}
\label{satz_Vert_Konv_1_2}
For $(i,j) \in \left\{ (1,2),(3,4)\right\}$ we have
$$
n^{2/3}\big( 2a - M_{N}^{i,j}\big)
\overset{\mathcal{D}}{\longrightarrow } \min_{k,l \in \N}\left\{Z_{1,k}^{i} + Z_{1,l}^{j} +\frac{a}{4}\left(Z_{2,k}^{i} + Z_{2,l}^{j}\right)^2 \right\}
$$
with two independent random elements
$(Z_{1,k}^{i},Z_{2,k}^{i})_{k \ge 1}  \sim \NA_{\infty}(\sigma_i,\tau_i)$ and $(Z_{1,k}^{j},Z_{2,k}^{j})_{k \ge 1}\sim \NA_{\infty}(\sigma_j,\tau_j).$
\end{env_prop}
\noindent
Now we are able to prove Theorem~\ref{satz_Vert_Konv_Hauptsatz} for $ \Xi_N$ (and thus, by \eqref{eq_Vert_Gleich_Zwei_diam} and \eqref{eq_Vert_Gleich_zwei_Xi}, also for $\Phi_n$).
\begin{proof} (of Theorem~\ref{satz_Vert_Konv_Hauptsatz})
\label{bew_Hauptsatz}
Let $
A_n := \big\{ \diam(\Xi_N) = \max\big{\{} M_{N}^{1,2} ,M_{N}^{1,3} ,M_{N}^{2,4} ,M_{N}^{3,4}\big{\}}\big\}
$
and fix $t >0$.
From  \eqref{eq_Stoch_Konv_diam_und_Max_von_4_max} we obtain
\begin{eqnarray*}
  & & \P\left(n^{2/3}\big( 2a - \diam({\Xi}_N)\big) \le t\right) \\
  & = &   \P\big(n^{2/3}\big( 2a - \diam({\Xi}_N)\big)\le t\big| A_n\big)\cdot \P(A_n) +  \P\big(n^{2/3}\big( 2a - \diam({\Xi}_N)\big)\le t\big| A_n^c\big)\cdot \P(A_n^c)\\
  & =&  \P\big(n^{2/3}\big( 2a - \max\left\{  M_N^{1,2},M_N^{1,3},M_N^{2,4},M_N^{3,4} \right\}\big)\le t\big)\big(1 + o(1)\big) + o(1).
\end{eqnarray*}
Since
$$
  n^{2/3}\big( 2a - \max\left\{  M_N^{1,2},M_N^{1,3},M_N^{2,4},M_N^{3,4} \right\}\big) \\
  =  \min_{(i,j) \in \{(1,2),(1,3),(2,4),(3,4)\}}\big\{n^{2/3}\left( 2a - M_N^{i,j}\right)\big\},
$$
the result follows from  Propositions~\ref{satz_Vert_Konv_1_3} to \ref{satz_Vert_Konv_1_2}.
\end{proof}


\section{Generalisations and open questions}
\label{sec_General_Open_Problems}
A very easy generalisation of our setting is given if we allow one boundary function to be equal to $0$ close to the corresponding pole.
Even the case $g_i(x_1) \equiv g_j(x_2) \equiv 0$ for $i \in \{1,4\}$, $j \in \{2,3\}$, $x_1$ close to $a$ and $x_2$ close to $-a$ is allowed. In these cases,
the minimum in \eqref{eq_Vert_Konv_Aussage} has to be taken over fewer random variables $S^{i,j}$, since certain combinations of quadrants do not contribute to
the maximum interpoint distance for $n$ large.
Because of the unique major axis, A6 can be weakened to a certain extent without changing the asymptotic behaviour of $M_n$. We omit the details.
Instead of A7, we can demand that there are constants $p_i >0$, $i \in \{1,2,3,4\}$, such that
$$
\sup \Big\{ \big|f(z) - p_i \big| : z \in E_i \cap U(a_i,\varepsilon)  \Big\} \overset{\varepsilon \rightarrow 0}{\longrightarrow }  0
$$
with $a_1 := a_4 := a$ and $a_2 := a_3 := -a$. Now it is possible that $p_1 \ne p_4$ or $p_2 \ne p_3$. In this situation, all results remain unchanged.
Another obvious generalisation is to consider $k \ge 2$ major axes of $E$ with no common endpoints. If the boundary of $E$ fulfils all assumptions at every endpoint
(after suitable rotations), we obtain a minimum over $k$ independent random variables $\min\left\{S_i^{1,2},S_i^{1,3}, S_i^{2,4},S_i^{3,4}\right\}$ as limit distribution (we omit the details).
A completely different setting is given if we weaken A7 by demanding $f(z)>0$ for $z$ close to the poles but $f(z) \to 0$ as $z$ tends to $(a,0)$ (resp. $(-a,0)$). Can $2a - \diam(\Phi_n)$ be scaled appropriately  (depending on the speed of the convergence to $0$) to obtain a limit distribution of $M_n$ also in this case? We leave this as an open problem.
\\

\noindent {\bf Acknowledgement:} The author wishes to thank Prof.~Dr.~Norbert Henze for bringing this problem to his attention and for helpful discussions.

\bibliographystyle{elsarticle-harv}
\bibliography{C:/Users/schrempp/Desktop/LitVerzeichnis}

\end{document}